\documentclass[12pt]{article}

\usepackage{latexsym,amsmath,amssymb,amsthm,amsfonts}
\usepackage{graphicx}
\usepackage{xcolor}
\usepackage{subfigure}
\usepackage{hyperref}

\title{\bfseries%
  Gradient descent in a generalised Bregman distance framework
}
\author{%
  Martin Benning$^1$, Marta M. Betcke$^2$, Matthias J. Ehrhardt$^1$,\\
  and Carola-Bibiane Sch\"{o}nlieb$^1$\\[12pt]
  \normalsize
  University of Cambridge, Dept. of Applied Mathematics \& Theoretical Physics,\\ \normalsize Wilberfore Road, Cambridge CB3 0WA, United Kingdom\\
  \normalsize
  \texttt{\{mb941, me404, cbs31\}@cam.ac.uk}\\[6pt]
  \normalsize
  University College London, Dept. of Computer Science, Gower Street,\\ \normalsize London WC1E 6BT, United Kingdom\\
  \normalsize
  \texttt{m.betcke@ucl.ac.uk}\\
}
\date{}

\DeclareMathOperator*{\argmin}{\arg \min}
\DeclareMathOperator{\dom}{dom}
\DeclareMathOperator{\uc}{\mathcal{U}}

\newcommand{\sbd}{D^{\text{symm}}_J}
\newcommand{\sbds}{D^{\text{symm}}_{J^*}}
\newcommand{\R}{\mathbb{R}}
\newcommand{\N}{\mathbb{N}}

\newtheorem{definition}{Definition}
\newtheorem{lemma}{Lemma}
\newtheorem{remark}{Remark}
\newtheorem{theorem}{Theorem}

\begin{document}

\maketitle
\pagestyle{empty}
\thispagestyle{empty}

%%
%% main text
\section{Introduction}
\label{sec:INTRODUCTION}
In this work we study a generalisation of classical gradient descent that has become known in the literature as the so-called linearised Bregman iteration \cite{yin2008bregman,yin2010analysis}, and -- as the key novelty of this publication -- apply it to minimise smooth but not necessarily convex objectives $E:\uc \rightarrow \R$ over a Banach space $\uc$. For this generalisation we want to consider proper, lower semi-continuous (l.s.c.), convex but not necessarily smooth functionals $J:\uc \rightarrow \R \cup \{ \infty \}$, and consider their generalised Bregman distances
\begin{align*}
D_J^p(u, v) = J(u) - J(v) - \langle p, u - v \rangle\, , %\label{eq:bregdis}
\end{align*}
for $u, v \in \uc$ and $p \in \partial J(v)$, where $\partial J(v)$ denotes the subdifferential of $J$. Note that in case $J$ is smooth we omit $p$ in the notation of the Bregman distance, as the subdifferential is single-valued in this case. We further assume that there exists a proper, l.s.c., convex and not necessary smooth functional $F:\uc \rightarrow \R \cup \{ \infty \}$ such that the functional $G := F - E$ is also convex. This will imply $D^{q - \nabla E(v)}_G(u, v) \geq 0$ for all $u, v \in \dom(G)$ and $q \in \partial F(v)$, since $q - \nabla E$ is the gradient of $G$. Hence, the convexity of $G$ yields the descent estimate
\begin{align}
E(u) \leq E(v) + \langle \nabla E(v), u - v \rangle + D_F^q(u, v) \, , \label{eq:descest}
\end{align}
for all $u, v \in \dom(F)$ and $q \in \partial F(v)$. We want to emphasise that in case of $F(u) = \frac{L}{2}\|u\|_{2}^2$ (for some constant $L > 0$) \eqref{eq:descest} reduces to the classical Lipschitz estimate; this generalisation has also been discovered in \cite{bauschkedescent} simultaneously to this work (without the generalisation of Bregman distances to non-smooth functionals, though).

\section{Linearised Bregman iteration applied to non-convex problems}
\label{sec:linbregan}
The linearised Bregman iteration that we are going to study in this work is defined as
\begin{subequations}\label{eq:linbregman}
\begin{align}
u^{k + 1} &= \argmin_{u \in \dom(J)} \left\{ \tau^k \langle u - u^k, \nabla E(u^k) \rangle + D_J^{p^k}(u, u^k) \right\} \, \text{,}\label{eq:linbregman1} \\
p^{k +1} &= p^k - \tau^k \nabla E(u^k) \, \text{,}\label{eq:linbregman2}
\end{align}
\end{subequations}
for $k \in \N$, some $u^0 \in \uc$ and $p^0 \in \partial J(u^0)$. Here $J:\mathcal{U} \rightarrow \R \cup \{ \infty \}$ is not only proper, l.s.c. and convex, but also chosen such that the overall functional in \eqref{eq:linbregman1} is coercive and strictly convex and thus, its minimiser well-defined and unique.

We want to highlight that this model has been studied for several scenarios in which $E$ is the convex functional $E(u) = \frac{1}{2}\| Ku - f \|_2^2$, for data $f$ and linear and bounded operators $K$ (cf. \cite{yin2008bregman,yin2010analysis}), for more general convex functionals $E$ and smooth $J$ in \cite{nemirovski1982problem,beck2003mirror}, as well as for the non-convex functional $E(u) = \frac{1}{2}\| K(u) - f \|_2^2$ for data $f$ and a smooth but non-linear operator $K$ in \cite{Bachmayr2009}. However, to our knowledge this is the first work that studies \eqref{eq:linbregman} for general smooth but not necessarily convex functionals $E$.

\section{A sufficient decrease property}
\label{sec:convprop}
We want to show that together with the descent estimate \eqref{eq:descest} we can guarantee a sufficient decrease property of the iterates \eqref{eq:linbregman} in terms of the symmetric Bregman distance. The symmetric Bregman distance $\sbd(u, v)$ (cf. \cite{burger2007error}) is simply defined as $\sbd(u, v) = D_J^q(u, v) + D_J^p(v, u) = \langle u - v, p - q \rangle$ for all $u, v \in \dom(J)$, $p \in \partial J(u)$ and $q \in \partial J(v)$.
\begin{lemma}[Sufficient decrease property]\label{lem:suffdecrease}
Let $E:\uc \rightarrow \R$ be a l.s.c. and smooth functional that is bounded from below and for which a proper, l.s.c. and convex functional $F:\uc \rightarrow \R \cup \{ \infty \}$ exists such that $G:= F - E$ is also convex. Further, let $J:\uc \rightarrow \R \cup \{ \infty \}$ be a proper, l.s.c. and convex functional such that \eqref{eq:linbregman1} is well defined and unique. Further we choose $\tau^k$ such that the estimate
\begin{align}
\rho D_J^{\text{symm}}(u^{k + 1}, u^k) \leq \frac{1}{\tau^k} D_J^{\text{symm}}(u^{k + 1}, u^k) -  D_{F}^{q^k}(u^{k + 1}, u^k)\label{eq:rhoestimate}
\end{align}
holds true, for all $k \in \N$, $q^k \in \partial F(u^k)$ and a fixed constant $0 < \rho < \infty$. Then the iterates of the linearised Bregman iteration \eqref{eq:linbregman} satisfy the descent estimate
\begin{align}
E(u^{k + 1}) + \rho D^{\text{symm}}_J(u^{k + 1}, u^k) \leq E(u^k) \, \text{.}\label{eq:sufficientdecrease}
\end{align}
In addition, we observe
\begin{align*}
\lim_{k \rightarrow \infty} D^{\text{symm}}_J(u^{k + 1}, u^k) = 0 \, \text{.}
\end{align*}
\begin{proof}
First of all, we easily see that update \eqref{eq:linbregman2}, i.e.
\begin{align*}
\tau^k \nabla E(u^k) + (p^{k + 1} - p^k ) = 0 \, ,
\end{align*}
is simply the optimality condition of \eqref{eq:linbregman1}, for $p^{k + 1} \in \partial J(u^{k + 1})$. Taking a dual product of \eqref{eq:linbregman2} with $u^{k + 1} - u^k$ yields
% \begin{align*}
%\tau^k \langle \nabla E(u^k), u^{k + 1} - u^k \rangle + D^{\text{symm}}_J(u^{k + 1}, u^k) = 0 \, \text{,}
%\end{align*}
%and hence we obtain
\begin{align}
\langle \nabla E(u^k), u^{k + 1} - u^k \rangle = -\frac{1}{\tau^k} D^{\text{symm}}_J(u^{k + 1}, u^k) \, \text{.}\label{eq:intermedeq}
\end{align}
Due to \eqref{eq:descest} we can further estimate
\begin{align*}
E(u^{k + 1}) \leq E(u^k) + \langle u^{k + 1} - u^k, \nabla E(u^k)\rangle + D_{F}^{q^k}(u^{k + 1}, u^k) \, \text{,}
\end{align*} 
for $q^k \in \partial F(u^k)$. Together with \eqref{eq:intermedeq} we therefore obtain
\begin{align*}
E(u^{k + 1}) + \frac{1}{\tau^k} D_J^{\text{symm}}(u^{k + 1}, u^k) -  D_{F}^{q^k}(u^{k + 1}, u^k) \leq E(u^k) \, \text{.}
\end{align*}
Using \eqref{eq:rhoestimate} then allows us to conclude
\begin{align*}
0 \leq \rho D^{\text{symm}}_J(u^{k + 1}, u^k) \leq E(u^k) - E(u^{k + 1}) \, \text{;}
\end{align*}
hence, summing up over all $N$ iterates and telescoping yields
\begin{align*}
\sum_{k = 0}^N \rho D^{\text{symm}}_J(u^{k + 1}, u^k) &\leq \sum_{k = 0}^N E(u^k) - E(u^{k + 1}) \, \text{,}\\
&= E(u^0) - E(u^{N + 1}) \, \text{,}\\
&\leq E(u^0) - \overline E < \infty \, \text{,}
\end{align*}
where $\overline E$ denotes the lower bound of $E$. Taking the limit $N \rightarrow \infty$ then implies
\begin{align*}
\sum_{k = 0}^\infty \rho D^{\text{symm}}_J(u^{k + 1}, u^k) < \infty \, \text{,}
\end{align*}
and thus, we have $\lim_{k \rightarrow \infty} D^{\text{symm}}_J(u^{k + 1}, u^k) = 0$ due to $\rho > 0$.
\end{proof}
\end{lemma}

\begin{remark}\label{rem:duallimit}
We want to emphasise that Lemma \ref{lem:suffdecrease} together with the duality $D_J^{\text{symm}}(u^{k + 1}, u^k) = D_{J^*}^{\text{symm}}(p^{k + 1}, p^k)$, for $p^{k+1} \in \partial J(u^{k+1})$ and $p^k \in \partial J(u^k)$, further implies
\begin{align*}
\lim_{k \rightarrow \infty} D^{\text{symm}}_{J^*}(p^{k + 1}, p^k) = 0 \, \text{,}%\label{eq:dualitylimit}
\end{align*}
and hence, a sufficient decrease property holds also for the dual iterates. Here $J^* :\uc^* \rightarrow \R \cup \{ \infty \}$ denotes the Fenchel conjugate of $J$, and $\uc^*$ is the dual space of $\uc$.
\end{remark}

\section{A global convergence statement}
For the following part we assume that both $J$ and $J^*$ are strongly convex w.r.t. the $\uc$- respectively the $\uc^*$-norm, i.e. there exist constants $\gamma > 0$ and $\delta > 0$ such that
\begin{align}
\gamma \| u - v \|_{\uc}^2 \leq \sbd(u, v) \qquad \text{and} \qquad \delta \| p - q \|_{\uc^*}^2 \leq \sbds(p, q)\label{eq:strongconvexity}
\end{align}
hold true for all $u, v \in \uc$ and $p, q \in \uc^*$. From Lemma \ref{lem:suffdecrease} and \eqref{eq:strongconvexity} we readily obtain
\begin{align}
\rho_1 \| u^{k+1} - u^k \|_{\uc}^2 \leq E(u^k) - E(u^{k+1}) \, ,\label{eq:sqnormbounds}
\end{align}
for $\rho_1 := \gamma/\rho$, which implies $\lim_{k \rightarrow \infty} \| u^{k + 1} - u^k \|_{\uc} = 0$. %and $\lim_{k \rightarrow \infty} \| p^{k + 1} - p^k \|_{\uc^*}$.

We follow \cite{bolte2014proximal} and establish a global convergence result by proving that the dual norm of the gradient is bounded by the iterates gap in addition to the already proven descent result \eqref{eq:sqnormbounds}. Together with a generalised Kurdyka-\L ojasiewicz property we will be able to prove a global convergence statement for \eqref{eq:linbregman}. 

Given \eqref{eq:strongconvexity}, we obtain the necessary iterates gap in the corresponding Banach space norm as an upper bound for the gradient in the dual Banach space norm, as follows.
\begin{lemma}[Gradient bound]
Let the same assumptions hold true as in Lemma \ref{lem:suffdecrease}, and let \eqref{eq:strongconvexity} be fulfilled. Then the iterates \eqref{eq:linbregman} satisfy 
\begin{align}
\| \nabla E(u^k) \|_{\uc^*} \leq  \rho_2 \| u^{k + 1} - u^k \|_{\uc} \, ,\label{eq:gradientbound}
\end{align}
for $\rho_2 := 1/(\delta \overline{\tau})$ and $\overline{\tau} := \inf_k \tau^k$.
\begin{proof}
As pointed out in Remark \ref{rem:duallimit}, we have the duality $D^{\text{symm}}_{J^*}(p^{k + 1}, p^k) = D^{\text{symm}}_{J}(u^{k + 1}, u^k)$ for the symmetric Bregman distances. Together with the duality estimate $\langle u, p \rangle \leq \| u \|_{\uc} \| p \|_{\uc^*}$ we therefore obtain
\begin{align*}
D^{\text{symm}}_{J^*}(p^{k + 1}, p^k) = \langle p^{k + 1} - p^k, u^{k + 1} - u^k \rangle \leq \| u^{k + 1} - u^k \|_{\uc} \| p^{k + 1} - p^k \|_{\uc^*} \, \text{.}
\end{align*}
Hence, using \eqref{eq:linbregman2} yields
\begin{align*}
\frac{D^{\text{symm}}_{J^*}(p^k - \tau^k \nabla E(u^k), p^k)}{\tau^k \| \nabla E(u^k) \|_{\uc^*}} \leq \| u^{k + 1} - u^k \|_{\uc} \, \text{.}
\end{align*}
Together with the $\delta$-strong convexity \eqref{eq:strongconvexity} and $\rho_2 := 1/( \delta \overline{\tau})$ we get \eqref{eq:gradientbound}.
\end{proof}
\end{lemma}
\begin{remark}
Note that we have to ensure $\overline{\tau} > 0$ in order to ensure $\rho_2 < \infty$. Due to \eqref{eq:rhoestimate} we can ensure this as long as $D_F^{q^k}(u^{k + 1}, u^k)$ is bounded from above for all $k \in \N$. 
\end{remark}

Before we can establish a global convergence result, we have to restrict the functionals $E$ to the following class of functionals satisfying a generalised Kurdyka-\L ojasiewicz property.
\begin{definition}[Generalised Kurdyka-\L ojasiewicz (KL) property]\label{def:kl}
We assume for $\eta > 0$ that $\varphi:[0, \eta [ \rightarrow \mathbb{R}_{> 0}$ is a function that is continuous at zero and satisfies $\varphi(0) = 0$, $\varphi \in C^1(]0, \eta [)$. Let further $E:\uc \rightarrow \R$ be a proper, l.s.c. and smooth functional.
\begin{enumerate}
\item The functional $E$ fulfils the (generalised) KL property at a point $\overline{u} \in \uc$ if there exists $\eta \in ]0, \infty]$, a neighbourhood $U$ of $\overline{u}$ and a function $\varphi$ satisfying the conditions above, such that for all
\begin{align*}
u \in U \cap \{ u \ | \ E(\overline{u}) < E(u) < E(\overline{u}) + \eta\}
\end{align*}
we observe
\begin{align}
\varphi^\prime(E(u) - E(\overline{u}))\| \nabla E(u) \|_{\uc^*} \geq 1 \, .\label{eq:modifiedkl}
\end{align}
\item If $E$ satisfies the (generalised) KL property for all arguments in $\uc$, $E$ is called a (generalised) KL functional.
\end{enumerate}
\end{definition}
%\begin{definition}[Generalised Kurdyka-\L ojasiewicz (KL) property]
%We assume for $\eta > 0$ that $\varphi:[0, \eta [ \rightarrow \mathbb{R}_{> 0}$ is a function that is continuous at zero and satisfies $\varphi(0) = 0$, $\varphi \in C^1(]0, \eta [)$, and $\varphi^\prime(s) > 0$ for all $s \in ]0, \eta [$. Then the proper and l.s.c. functional $E$ is said to satisfy the generalised KL property at the critical point $\hat{E} \in \uc$ - this means $\nabla E(\hat{u}) = 0$ - if
%\begin{align}
%\varphi^\prime\left(E(u) - E(\hat{u})\right)\| \nabla E(u) \|_{\uc^*} \geq 1 \, \text{.}\label{eq:modifiedkl}
%\end{align}
%holds true for all $u \in \uc$.
%\end{definition}
Together with the previous results the generalised KL condition \eqref{eq:modifiedkl} allows to establish the following global convergence result.
\begin{theorem}[Global convergence]
Let the Banach space $\uc$ be the dual of a separable normed space. Suppose that $E$ is coercive, sequentially weak$^\ast$-continuous and a KL function in the sense of Definition \ref{def:kl}. Then the sequences $\{ u^k \}_{k \in \N}$ and $\{ p^k \}_{k \in \N}$ generated by \eqref{eq:linbregman} each have a strongly convergent subsequence with limits $\hat u$ and $\hat p$, with $\nabla E(\hat u) = 0$ and $\hat{p} \in \partial J(\hat{u})$. If $\text{dim}(\uc) < \infty$, then the convergence holds true for the entire sequences.
%\begin{enumerate}
%\item have finite length, i.e.
%\begin{align}
%\sum_{k = 0}^\infty \| u^{k + 1} - u^k \|_{\uc} < \infty \qquad \text{and} \qquad \sum_{k = 0}^\infty \| p^{k + 1} - p^k \|_{\uc^*} < \infty \, \text{.}\label{eq:finitelength}
%\end{align}
%\item converge to a critical point $\hat{u}$ of $E$ with $\hat{p} \in \partial J(\hat{u})$.
%\end{enumerate}
\begin{proof}
The proof utilises \eqref{eq:sufficientdecrease}, \eqref{eq:gradientbound} and \eqref{eq:modifiedkl} to derive the statement. Due to page restrictions, the full length proof will be published separately in an extended version of this manuscript.
\end{proof}
\end{theorem}

\begin{figure}[hbt]
  \centering
    \subfigure[Ground truth]{\includegraphics[width=0.24\textwidth]{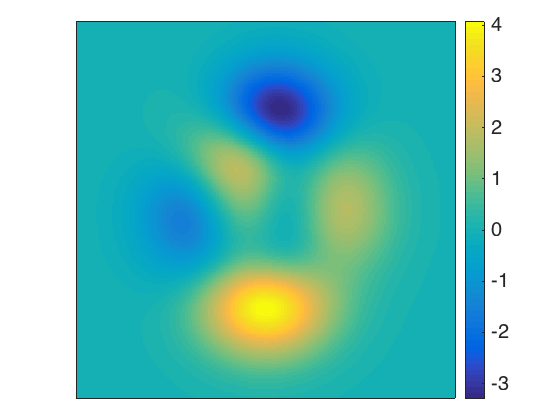}\label{subfig:gt}}
    %\subfigure[Wrapped signal]{\includegraphics[width=0.32\textwidth]{Images/wrapped.png}\label{subfig:wrapped}}\\
    \subfigure[Gradient descent]{\includegraphics[width=0.24\textwidth]{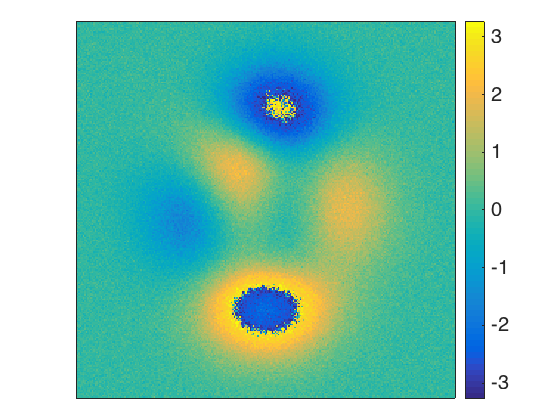}\label{subfig:graddescent}}
    \subfigure[\mbox{$R(u) = \frac{1}{2} \| \nabla u \|_{L^2}^2$}]{\includegraphics[width=0.24\textwidth]{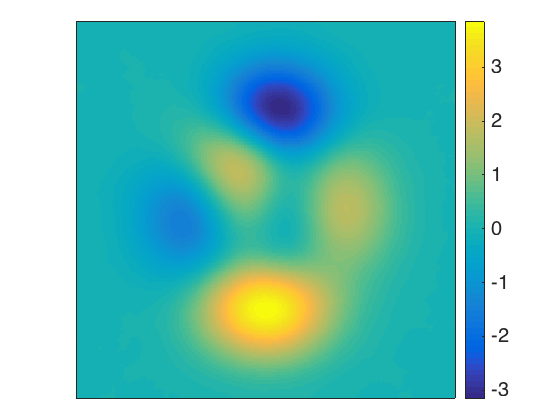}\label{subfig:h1descent}}
    \subfigure[$R(u) = \| Cu \|_{\ell^1}$]{\includegraphics[width=0.24\textwidth]{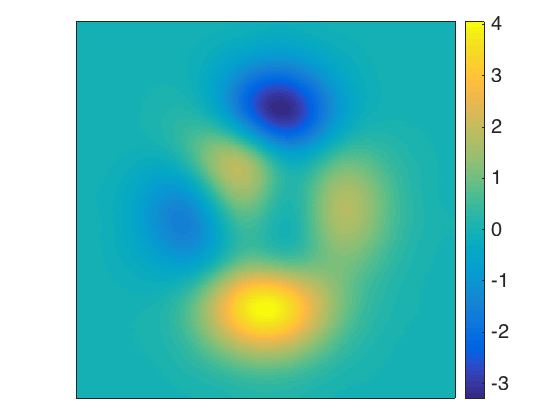}\label{subfig:dctdescent}}
    \caption[Phase unwrapping]{A phase unwrapping example. Figure \ref{subfig:gt} shows the unknown, noise-free, ground truth signal. Figure \ref{subfig:graddescent} shows the result of classical gradient descent computation. Figure \ref{subfig:h1descent} visualises the solution of model 2.) with $\alpha = 1000$. Figure \ref{subfig:dctdescent} shows the solution of model 3.) with $\alpha = 50$. All reconstructions have been computed from zero initialisations and were stopped according to the same discrepancy principle.}
    \label{fig:phasetoyexm}
\end{figure}

\section{Phase unwrapping as a toy example}
\label{sec:results}
We want to conclude this paper with a numerical toy example for which we consider to minimise $E(u) := \frac{1}{2} \| K(u) - f \|_{L^2(\Omega; \R^2)}^2$ for $K(u) = ( \cos(u), \sin(u) )^T$, and choose $F(u) = \frac{L}{2}\| u \|_{L^2(\Omega)}^2$ with $L = 1$. We will minimise $E$ via \eqref{eq:linbregman} with $J(u) := \frac{1}{2}\| u \|_{L^2(\Omega)}^2 + \alpha R(u)$, for a positive scalar $\alpha > 0$ and three different choices of $R$: 1.) $R(u) = 0$, 2.) $R(u) = \frac{1}{2}\| \nabla u \|_{L^2(\Omega; \R^2)}^2$, and 3.) $R(u) = \| Cu \|_{\ell^1}$, where $C$ denotes the two-dimensional discrete Cosine transform. The first case simply corresponds to classical gradient descent, case 2.) is gradient descent in a Hilbert space metric and 3.) corresponds to gradient descent in a non-smooth Bregman distance setting that does not correspond to a metric. Note that the question, whether $E$ and $J$ satisfy all conditions that are necessary for global convergence, will be omitted due to the page limit, but addressed in an extended version of this manuscript in the future. We do want to mention, though, that it is easy to see that $J$ in 3.) does not meet the requirement \eqref{eq:sqnormbounds}; this, however, can be corrected via a smoothing of the $\ell^1$-norm, for instance via a Huberised $\ell^1$-norm.

In order to consider numerical examples, we discretise the above scenarios in a straight forward fashion. Input data $f$ is created by applying the non-linear operator $K$ to a multiple of the built-in MATLAB\copyright{} signal 'peaks' (see Figure \ref{subfig:gt}) and additive normal distributed noise with mean zero and standard deviation $\sigma = 0.15$. Due to noise in the data, the iteration \eqref{eq:linbregman} is stopped as soon as $E(u^k) \leq \sigma^2 m/2$ is satisfied. Here $m$ denotes the number of discrete samples. Reconstruction results for zero initialisations and the choice $\tau^k = 1.5$ for all $k \in \N$ can be found in Figure \ref{subfig:graddescent}, \ref{subfig:h1descent} and \ref{subfig:dctdescent}. We want to emphasise that this example is just a toy example to demonstrate the impact of different choices of $J$; there are certainly much better unwrapping strategies, particularly for the unwrapping of smooth signals.\\ 

\noindent\textbf{Code statement}: The corresponding MATLAB\copyright{} code can be downloaded at\\ \url{https://doi.org/10.17863/CAM.6714}.

%\section{Preliminary numerical results}
%\label{sec:results}

\section{Conclusions \& Outlook}
We have presented a short convergence analysis of the linearised Bregman iteration for the minimisation of general smooth but non-convex functionals. We have proven a sufficient decrease property, and confirmed that the dual norm of the gradient is bounded by the primal iterates under additional strong convexity assumptions of the convex functional that builds the basis for the Bregman iteration. Under a generalised KL condition, we have stated a global convergence result that we are going to refine in detail in a future release. We have concluded with a numerical toy example of phase unwrapping for three different Bregman distances. 
In a future work we are going to analyse the linearised Bregman iteration and its convergence behaviour in more detail and in a more generalised setting, and are going to investigate different Bregman distance choices as well as different numerical applications.

%%
%% acknowledgment
\section*{Acknowledgment}
MB acknowledges support from the Leverhulme Trust Early Career project 'Learning from mistakes: a supervised feedback-loop for imaging applications' and the Isaac Newton Trust. MMB acknowledges support from the Engineering and Physical Sciences Research Council (EPSRC) 'EP/K009745/1'. MJE and CBS acknowledge support from the Leverhulme Trust project 'Breaking the non-convexity barrier'. CBS further acknowledges support from EPSRC grant 'EP/M00483X/1', EPSRC centre 'EP/N014588/1' and the Cantab Capital Institute for the Mathematics of Information. All authors acknowledge support from CHiPS (Horizon 2020 RISE project grant) that made this contribution possible. 

%%
%% bibliography
%\begin{thebibliography}{9}
%
%\bibitem{AA1-2012a}
%\newblock Author, A.,
%\newblock {\em Book Title\/},
%\newblock Publisher, Place, 2012.
%
%\bibitem{AA2-2012a}
%\newblock Author, A. and Author, B.,
%\newblock {\em Post-conference paper title\/},
%\newblock in Editor, A. and Editor, B., eds.,
%          Proceedings of the Conference Name, Place,
%\newblock Vol.6, No.1 (2002), pp.1--4.
%
%\bibitem{AA3-2012a}
%\newblock Author, A., Author, B., and Author, C.,
%\newblock {\em Paper title\/},
%\newblock Journal Name, 6 (2002), pp.1--44.
%
%\end{thebibliography}

\bibliography{bib}
\bibliographystyle{plain}

\end{document}